\documentclass[10pt,amssymb]{article}
\usepackage{amsfonts,amsmath,latexsym,verbatim,amscd}

\title{Maskit combinations of \\ Poincar\'e-Einstein metrics}

\author{Rafe Mazzeo\thanks{Email: mazzeo@math.stanford.edu. 
Supported by the NSF under Grant DMS-0204730}\\ Stanford University \and  
Frank Pacard \thanks{Email: pacard@univ-paris12.fr}\\ 
Universit\'e Paris XII}

\date{}

\begin{document}
\maketitle

\newcommand{\CC}{\mathbb C}
\newcommand{\HH}{\mathbb H}
\newcommand{\RR}{\mathbb R}
\newcommand{\VV}{{\mathbb V}}
\newcommand{\del}{\partial}
\newcommand{\e}{\varepsilon}
\newcommand{\al}{\alpha}
\newcommand{\olg}{\overline{g}}
\newcommand{\Ric}{\mathrm{Ric}}
\newcommand{\calA}{{\mathcal A}}
\newcommand{\calC}{{\mathcal C}}
\newcommand{\calD}{{\mathcal D}}
\newcommand{\calE}{{\mathcal E}}
\newcommand{\calF}{{\mathcal F}}
\newcommand{\calG}{{\mathcal G}}
\newcommand{\calL}{{\mathcal L}}
\newcommand{\calM}{{\mathcal M}}
\newcommand{\calO}{{\mathcal O}}
\newcommand{\calP}{{\mathcal P}}
\newcommand{\calS}{{\mathcal S}}
\newcommand{\calT}{{\mathcal T}}
\newcommand{\calU}{{\mathcal U}}
\newcommand{\calV}{{\mathcal V}}
\newcommand{\Conf}{{\mathrm{Conf}}}
\newcommand{\tr}{{\mathrm{tr}}\,}
\newcommand{\tg}{\tilde{g}}
\newcommand{\thh}{\tilde{h}}
\newcommand{\tw}{\tilde{w}}
\newcommand{\tx}{\tilde{x}}
\newcommand{\ty}{\tilde{y}}
\newcommand{\tir}{\tilde{r}}
\newcommand{\tphi}{\tilde{\phi}}
\newcommand{\tom}{\tilde{\omega}}
\newcommand{\Mbar}{\overline{M}}
\newcommand{\dMbar}{\del\Mbar}
\newcommand{\frakc}{{\mathfrak c}}
\newcommand{\PE}{Poincar\'e-Einstein\ }

\newtheorem{theorem}{Theorem}
\newtheorem{proposition}{Proposition}
\newtheorem{corollary}{Corollary}
\newtheorem{lemma}{Lemma}
\newtheorem{definition}{Definition}
\newtheorem{remark}{Remark}

\begin{abstract}
We establish a boundary connected sum theorem for asymptotically hyperbolic Einstein 
metrics; this requires no nondegeneracy hypothesis. We also show that if the two metrics have 
scalar positive conformal infinities, then the same is true for this boundary join.  
\end{abstract}

\section{Introduction}

Let $M^{n+1}$ be a compact manifold with boundary. A Riemannian metric $g$ on the {\it interior} 
of $M$ is said to be conformally compact if $g = \rho^{-2}\olg$, where $\olg$ is nondegenerate 
(with some  specified regularity) up to the boundary and $\rho$ is a defining function for $\del M$ 
(i.e.\ $\rho^{-1}(0) = \del M$ and  $d\rho \neq 0$ there). Any such $g$ is complete, and if 
$|d\rho|_{\olg} = 1$ on $\del M$, then its sectional curvatures tend  to $-1$ as $\rho \to 0$.
These metrics generalize the Poincar\'e model  of hyperbolic space $\HH^{n+1}$, and accordingly,
we call $(M,g)$  Poincar\'e-Einstein (or simply PE) if $g$ is both conformally compact 
and Einstein. We always normalize so that $\Ric^g = -ng$.

Poincar\'e-Einstein metrics, which are also known as asymptoticaly hyperbolic Einstein (AHE) or conformally
compact Einstein metrics, were introduced originally by Fefferman  and Graham \cite{FG} as a tool in conformal geometry. More recent interest in them has been generated by their role in the AdS/CFT correspondence in string theory, and this has stimulated much interesting work in this area. A number of explicit examples are known, starting from the most elementary one, the hyperbolic space $\HH^{n+1}$ and its convex cocompact quotients $\HH^{n+1}/\Gamma$, but also including the hyperbolic analogue of the Schwarzschild metric and, when $\dim M = 4$, the Taub-BOLT metrics on disk bundles over Riemann surfaces. Many of these are catalogued in \cite{A3}. 

More general \PE metrics can be obtained by perturbing known examples,  as established first in the work of Graham and Lee \cite{GL}, and much  later, in more general circumstances, by Biquard \cite{Bi} and Lee  \cite{L}. A substantially more comprehensive theory, especially in  four dimensions, has been obtained recently by Anderson \cite{A1},  \cite{A2}, \cite{A3}. He shows first that when $\dim M$ is arbitrary,  the moduli space $\calE = \calE(M)$ of \PE metrics on $M$ is either  empty or else a smooth infinite dimensional Banach manifold (with  respect to a suitable Sobolev or H\"older completion). Unlike the  situation when $M$ is closed without boundary, this deformation theory  is always unobstructed and $\calE(M)$ has no singularities.
In  addition, when $\dim M = 4$ it possible to set up a degree theory to  obtain more global existence
results. We explain this more carefully in the next section.

The basic goal of all of these papers is to solve an asymptotic  boundary problem. More specifically, there is a map $\frakc$  which associates to any conformally compact metric its {\it conformal infinity} $\frakc(g)$, which is by definition the conformal class of the restriction of $\olg = \rho^2 g$ to $\del M$.    A preliminary conjecture is that the restriction of $\frakc$ to  $\calE$ is a bijection, or at least a surjection; in other words,  every conformal class on $\del M$ is the conformal infinity of at least one \PE  metric. Subject to a certain nondegeneracy condition, this is true locally, i.e.\  if $g$ satisfies this nondegeneracy condition, then all conformal classes  near to $\frakc(g)$ are in the image of $\frakc$. However, explicit metrics are known which do not satisfy this condition \cite{A3};  for example, the conformal class of the product metric $S^{n-1} \times S^1(L)$ is not in the image of  $\frakc$ for $M = B^n \times S^1(L)$ when
the length $L$ of the $S^1$  factor is sufficiently large. Anderson shows \cite{A3} that 
$\frakc: \calE(M) \to \Conf\,(\del M)$ is always Fredholm of index zero. 

Our goal in this paper is to construct a wider class of \PE metrics  using a method which generalizes the boundary connected sum procedure  in hyperbolic geometry; this hyperbolic construction is part of the  Maskit combination theorems. More specifically, suppose that  $(M_j,g_j)$ are two \PE metrics. Fix points $p_j \in \del M_j$ and  excise small half-balls $B^{n+1}_+(p_j)$. The boundary connected sum  $M_1 \#_b M_2$ is obtained by identifying the hemispherical  portions of these boundaries. In the following we let $B_{j,+}$ denote any such half-ball in $M_j$ centered at $p_j$.

\begin{theorem} 
The manifold $M = M_1 \#_b M_2$ carries a family of  \PE  metrics $g_\e$ with the following two properties: 
\begin{itemize}
\item the restriction of $g_\e$ to $M_j - B_{j,+}$ converges to $g_j$;
\item the restriction of $ \rho^2 \, g_\e$ to $\del M_j - (B_{j,+} 
\cap \del M_j)$ converges to $\rho^2 \, g_j$.
\end{itemize}
The convergence in either case is polynomial in a geometrically natural parameter $\e$.  
\end{theorem}
In the case where the metric $g_i$ are nondegenerate (see  Definition~\ref{de:ndg}), the second statement can be improved since,  for the metrics we construct, the restriction of $\rho^2 \, g_\e$ to  $\del M_j - (B_{j,+} \cap \del M_j)$ is identically equal to  $\rho^2 \, g_j$. In particular, this implies that $\frakc (g_\e)$ is  equal to $\frakc (g_i)$ on $\del M_j - (B_{j,+} \cap \del M_j)$.

An important theme in this theory is that \PE metrics which have conformal infinities of positive (or at least nonnegative) Yamabe type are geometrically more stable, and the existence theory (in dimension four) is certainly more robust in this case. Recall that  a conformal class is said to be positive (nonnegative) if it  contains a metric with positive (nonnegative) scalar curvature.  As an example of this stability, the main step in Anderson's development of a ${\mathbb Z}$-valued degree theory for $\frakc$ in dimension $4$ is the  properness for the restriction of this map to the preimage of $\Conf^+(\del M)$,  the space of positive conformal classes.  Another example, in general dimensions,  is the well-known result of Witten and Yau \cite{WY}, cf.\ also \cite{CG}, that if $M$ carries a \PE metric with positive conformal infinity, then $H_n(M,\del M) = 0$ (and so, in particular, $\del M$ is connected). A primary motivation for our Theorem~1 is to construct many more \PE metrics with positive conformal infinities. 

\begin{corollary} 
Suppose that $(M_j,g_j)$, are \PE and furthermore,  that $\frakc(g_j) \in \Conf^+(\del M_j)$, $j=1,2$. Then  $\frakc(g_\e) \in \Conf^+(M_1 \#_b M_2)$ when $\e$ is small. In  particular, when $\dim M = 4$, every $3$-manifold which is a finite  connected sum of quotients of $S^3$ and $S^2 \times S^1$ arises as the  boundary of a \PE metric with positive conformal infinity. 
\end{corollary}

A well-known theorem of Schoen and Yau \cite{SY}, cf.\ also  \cite{GrL}, states that an arbitrary $3$-manifold which carries a  metric of positive  scalar curvature is a connected sum of quotients  of $S^2 \times S^1$ and manifolds with finite fundamental group.  Contingent on two (big!!) conjectures (the Poincar\'e and  the spherical space-form conjecture), these latter summands are all  lens spaces, and hence are known 
\cite{CS} to have explicit \PE fillings.  Thus, modulo these conjectures, or more simply, if one could show that  positive conformal classes on quotients of exotic homotopy $3$-spheres either  have \PE fillings, or else do not exist, then our Corollary~1 would imply  that every $3$-manifold $Y$ which admits a  metric of positive scalar curvature has a \PE filling with positive  conformal infinity. In any case, combining the results here with  Anderson's, we obtain a large class of new examples of \PE manifolds ( in arbitrary dimension) with positive conformal infinity. 

One can also ask whether other surgery constructions are  possible in  this category. The most direct generalization would be to join the  manifolds $M_j$ along a common submanifold $\Sigma \hookrightarrow \del M_j$  (with isomorphic normal bundles $N_j(\Sigma) \subset T_\Sigma \del M_j$).  Unfortunately, our construction  does not go through in any direct  way, and may even fail. This is not to say that the manifold  $M_1 \#_\Sigma M_2$ does not carry any \PE metrics, but if these  exist, they seem to be distant from the initial metrics $g_j$, even well away from a neighbourhood of $\Sigma$. 
We discuss this further at the end of this paper.

\S 2 contains a review of some details of the geometry of \PE metrics  and of the analytic methods used in their perturbation theory. This is  followed in three subsequent  sections by the construction of  approximate solutions, the linear estimates and the proof of the main  theorem.  In \S 6 we prove Corollary~1, and in \S 7 we discuss the plausibility of more general gluing constructions.  

\section{Preliminaries}

We now review in more detail some of the geometric and analysis required in our study of \PE metrics. 

\subsection{Geometry of \PE metrics}

Suppose that $g$ is a conformally compact metric on $M$, so that it  can be written  $g = \rho^{-2}\olg$ for some defining function $\rho$.  We shall always suppose that both $\rho$ and $\olg$ are (at least)   $\calC^{2,\alpha}(\overline{M})$. The precise regularity is not so  important in this paper, but we shall address this issue more  carefully below.

\subsubsection{Boundary normal coordinates}

It will be very convenient to have something like Fermi coordinates  around $\del M$. Following 
Graham and Lee \cite{GL}, we have
\begin{lemma}[\cite{GL}]
Assume that $g$ is a conformally compact metric so that $g = \rho^{-2}\, \olg$ for some defining 
function $\rho$, and $h_0$ is a representative of  $\frakc (g)$, the conformal infinity of $g$. 
Further assume that  $|d\rho|^2_{\bar g}=1$ at $\rho=0$. Then there exists a new defining function 
$x$, in terms of which the metric $g$ can be written as 
\begin{equation}
g = \frac{dx^2 + h}{x^2}
\label{eq:fc-cc}
\end{equation}
in some neighborhood of $\del M$. Here $h = h(x)$ is a family of Riemannian metrics on 
$\del M$ depending parametrically on $x$ with $h(0)=h_0$. 
\label{le:1}
\end{lemma}
Since this result is crucial to our construction, we briefly recall its proof now.  If $\rho$ is an initial defining function, then we look for a new one of the form $x = e^{u}\, \rho$. The metric $g$ will be in the correct form provided $|dx|_{x^2 g} \equiv 1$ near $\del M$, and this is equivalent to the equation for the function $u$
\begin{equation}
2 \, < d\rho, du>_{\bar g} \ + \ \rho \, |du|^2_{\bar g} \ = \ \frac{1 - |d \rho|^2_{\bar g}}{\rho}.
\label{eq:eqdf}
\end{equation}
Since $d\rho \neq 0$ and $\rho=0$ on $\del M$, equation (\ref{eq:eqdf}) is  noncharacteristic with respect to the boundary, and hence can be solved locally with any boundary condition $u=u_0$ when $\rho=0$. We take $u_0$ so that $e^{-2u_0}\, \olg|_{\del M} = h_0$. In terms of $x = e^{u} \, \rho$ and any choice of coordinates on $\del M$,   the metric $g$ is given by (\ref{eq:fc-cc}). 

\medskip

Recall from \cite{Be} that if two metrics are conformally related, say $g = e^{2f} \olg$, then the Ricci tensors of $g$ and $\olg$ are related by
\begin{equation}
\Ric^g = \Ric^{\olg} - (n-1) \, (\nabla^{\olg} df - df \circ df) - (\Delta^{\olg} f + (n-1) \, |df|^2_{\olg} ) \, \olg.
\label{eq:ricg}
\end{equation}
Applying this with $f = -\log x$, and the metric $\olg = dx^2 + h$  given in Lemma~\ref{le:1}, we obtain the expansion
\begin{equation}
\Ric^{g} + n \, g = \frac{1}{2x}\left( (n-1)h_{1} +  (\tr^{h_{0}}(h_{1})) \, (dx^2+ h_0)\right) + {\cal O} (1).
\label{eq:esterr}
\end{equation}
where we have expanded the metric  $h = h_0 + x \, h_1 + {\cal O}(x^2)$. As a consequence we have the~:
\begin{lemma}[\cite{GL}]
Under the assumptions and notations of Lemma~\ref{le:1}, if we further  assume that $g$ is a Poincar\'e-Einstein metric, then the family of  Riemannian metrics $h(x)$ on $\del M$ can be expanded as  $h(x) = h_0 + x^2 \, \tilde h(x)$, where $\tilde h(x)$ is a family of  symmetric $2$-tensors on $\del M$ which depend parametrically on $x$ .
\label{le:2}
\end{lemma}

In other words, the Einstein condition implies that $h'(0) = 0$ (or equivalently, that $\del M$ is 
totally geodesic for the metric  $\olg = x^2 g$), and hence $h = h_0 + x^2 \, h_2 + {\it o}(x^2)$.   
Assuming that $\rho$ and $\olg$ are polyhomogeneous, Fefferman and  Graham \cite{FG} produce a 
complete {\it formal} expansion for $h$; this is justified in four dimensions by the regularity 
result in  \cite{A2}, and in general dimensions in the forthcoming paper \cite{GM}. 

\subsubsection{Gauge choice and the Einstein equation}

The equation satisfied by $g$ is $\Ric^g + ng = 0$. It is  well-known that this equation is not elliptic because of the underlying diffeomorphism invariance, and so one must choose some gauge condition. The best choice is the one adopted by Biquard \cite{Bi}, and later Anderson, called the Bianchi gauge. If $\tg = g+k$ is any metric near to $g$, then we define
\begin{equation}
B^g(\tg) = \delta^g \tg + \frac12 \, d \, \tr^g \, \tg = \delta^g k + \frac12 d\, tr^g \, k
\label{eq:bg}
\end{equation}
since $B^g(g)=0$. Thus $B^g$ is a map from symmetric $2$-tensors to $1$-forms. The  system
\[
\Ric^{\tg} + n \tg = 0, \qquad B^g(\tg) = 0 ,
\]
is elliptic in the sense of Agmon-Douglis-Nirenberg, but it is more convenient to work with the single equation 
\begin{equation}
N_g(k) := \Ric^{g+k} + n(g+k) + (\delta^{g+k})^*B^g(k) = 0,
\label{eq:main}
\end{equation}
where the symmetric $2$-tensor $k$ is assumed to be small enough so that $g+k$ is a metric on $M$.
\begin{proposition}[\cite{Bi}]
Suppose that $\Ric^{g+k} < 0$ and $|B^g(g+k)| \to 0$ at $\del M$; then any solution of $N_g(k) = 0$ corresponds to an Einstein metric $g+k$ which is in the Bianchi gauge relative to  $g$.
\end{proposition}

The linearization of $N_g$ is very simple in this gauge:
\begin{equation}
L_g \kappa := \left. 2 DN_g \right|_0(\kappa) = (\nabla^g)^* \nabla^g  \ \kappa - 2 \stackrel{\circ}{R^g} \kappa \ +  \Ric^g \circ \kappa + \kappa \circ \Ric^g + 2 \, n \, \kappa. 
\label{eq:lin}
\end{equation}
Here 
\[
(\stackrel{\circ}{R^g} \kappa)_{ij} = R_{ipjq}\, \kappa^{pq}, \qquad \Ric^g \circ \kappa = \Ric_i^{\ p}\, \kappa_{pj}, \qquad \kappa \circ \Ric^g = \kappa_i^{\ p} \,\Ric_{pj},
\]
and all curvatures are computed relative to $g$.  Note that when $g$ is Einstein, $\Ric^g = -ng$, and hence
\[
L_g = (\nabla^g)^* \nabla^g - 2 \stackrel{\circ}{R^g} 
\] 
then.

This operator $L_g$ is not uniformly elliptic on $M$; rather it has the structure of a uniformly 
degenerate operator, as studied in detail in \cite{Ma1}, \cite{Ma2}. We shall require some of the main results of the theory of uniformly degenerate operators, and so we now digress briefly to explain this setup. These general results will either be stated for, or immediately specialized to, the operator $L_g$. 

\subsection{Uniformly degenerate operators}

Choose coordinates $w = (x,y): = (x,y_1, \ldots, y_n)$, where $x = w_0 \geq 0$ is a boundary defining function, near a point $p \in \del M$. A second order operator $L$ is said to be uniformly degenerate if it can be expressed in the form 
\begin{equation}
L = \sum_{j+|\alpha| \leq 2} a_{j,\alpha}(x, y)(x \, \del_{x})^j (x \, \del_{y})^\alpha,
\label{eq:unifdeg}
\end{equation}
where the (scalar or matrix-valued) coefficients are bounded. We usually assume that these coefficients are smooth on $\overline{M}$, but it is easy to extend most of the main conclusions of this theory when they are only polyhomogeneous, or of some finite regularity. Operators of this type arise naturally in geometry, and in particular all of the natural geometric operators associated to a conformally compact metric are uniformly degenerate. 

A typical example of such an operator is given by the Laplace-Beltrami operator on hyperbolic space $({\mathbb H}^{n+1}, g_0)$. Taking coordinates $(x,y) \in (0, \infty)\times {\mathbb R}^n$ in the upper half space model, we have 
\[
g_0 = \frac{dx^2+ dy^2}{x^2},
\]
and we obtain
\begin{equation}
\Delta_{g_0} =  x^2 \, \del_x^2 + (1-n) \, x \, \del_x + x^2 \, \Delta_y.
\label{eq:lhs}
\end{equation}

\subsubsection{Ellipticity and model operators}

There is a well-defined symbol in this context, 
\[
\sigma(L)(x,y;\xi,\eta) : = \sum_{j+|\alpha|=2} a_{j,\alpha}(x , y)\, \xi^j \eta^\alpha \neq 0 \qquad \mbox{when}\quad (\xi,\eta) \neq 0, 
\]
and we say that $L$ is elliptic in the uniformly degenerate calculus provided $\sigma(L)(x,y;\xi,\eta)$ is invertible when $(\xi,\eta) \neq 0$. 

Ellipticity alone is not enough to guarantee that $L$ is Fredholm between appropriate function spaces; 
one must also require that certain simpler operators which model $L$ near the boundary be invertible. 
There are two of these model operators: 
\begin{itemize}
\item The {\it normal operator} of $L$ is defined by 
\[
N(L) : = \sum_{j+|\alpha|\leq 2} a_{j,\alpha}(0, y) (s \del_{s})^j (s \del_{u})^\alpha, \qquad (s,u) \in \RR^+ \times \RR^n ;
\]
here $y\in \del M$ enters only parametrically and the operator acts on  functions on a half-space $\RR^+ \times \RR^n $, which is naturally identified with the  inward-pointing half tangent space $T^+_{(0, y)}M$. 
\item The {\it indicial operator} of $L$ is defined by
\[
I(L) : = \sum_{j \leq 2} a_{j,0}(0,y) (s\del_s)^j.
\]
\end{itemize}

For example, the normal and indicial operators associated to the (scalar) Laplace-Beltrami operator on 
hyperbolic space $({\mathbb H}^{n+1}, g_0)$ are given by
\[
N(\Delta_{g_0}) = s^2 \, \del_s^2 + (1-n) \, s \, \del_s + s^2 \Delta_u, \qquad
I(\Delta_{g_0}) = s^2 \, \del_s^2 + (1-n) \, s \, \del_s.
\]

The normal operator can be regarded as $L$ with its coefficients frozen (in an appropriate sense)
at the boundary, so the following result is not surprising. 

\begin{proposition}
If $g =  \rho^{-2}\, \olg$ is a smooth conformally compact metric such that $|d\rho|_{\olg}^2 = 1$ on $\del M$, then its Laplace-Beltrami operator $\Delta_g$ and the linearization $L_g$ of the gauged Einstein equation are both elliptic uniformly degenerate operators. Furthermore, their normal operators are given by
\begin{equation}
N(\Delta_g)  = \Delta_{g_0} \qquad N(L_g) = (\nabla^{g_0})^* \nabla^{g_0} - 2 \, \stackrel{\circ}{R^{g_0}},
\label{eq:nolap}
\end{equation}
respectively, where $g_0$ is the standard hyperbolic metric on $\HH^{n+1}$.
\end{proposition}

The indicial operator is a much more primitive model, but it captures  some fundamental invariants associated to $L$. 
\begin{definition} 
The number $\zeta \in \CC$ is said to be an indicial root of $L$ if 
\[
L(x^\zeta v(y)) = \calO(x^{\zeta+1}) \qquad \mbox{for any} \quad v \in \calC^\infty(\del M).
\]
\end{definition}
We may replace $L$ by $I(L)$ (and $x$ by $s$) here, since the higher order terms
which appear in $L$ but not in the indicial operator are irrelevant for this calculation. Thus
\[
I_L(\zeta)s^\zeta : = I(L)(s^\zeta) = \left(\sum_{j \leq 2} a_{j,0}(0,y)\zeta^j\right) s^\zeta ,
\]
and hence $\zeta$ is an indicial root of $L$ if and only if it is a root of the (matrix-valued) 
polynomial $\zeta \to I_L(\zeta)$. 

The operator $L$ acts naturally on weighted H\"older (and Sobolev)  spaces, and the indicial roots of $L$ determine the weights for which  these mappings do not have closed range.

The calculation of the indicial roots for the linearized Einstein operator $L_g$ is carried out efficiently in \cite{GL}, see also \cite{L}. There are several pairs  of indicial roots, corresponding to the action of $L_g$ on pure-trace and trace-free symmetric $2$-tensors; the latter space in turn decomposes into three 
irreducible summands, corresponding to the normal and tangential components of the $2$-tensors. 

Before giving the values of these indicial roots, we must fix a basis of sections for the space of 
symmetric $2$-tensors. There are two natural choices: the standard one, consisting of all symmetric 
products $dw_i\, dw_j$, and another, consisting of all symmetric  products $\frac{dw_i}{x}\frac{dw_j}{x}$. 
We use this latter choice since it is geometrically more natural -- the $1$-forms 
$dw_i/x$ are of length bounded away from infinity and zero with respect to any conformally compact 
metric $g$ -- and so we write any symmetric $2$-tensor $k$ as 
\[
k =  \sum_{i,j=0}^n k_{ij} \, \frac{dw_i}{x} \, \frac{dw_j}{x}.
\]
This is in accord with the notation and terminology of \cite{Ma1}, where the role of this
normalization is emphasized and exploited consistently. In fact, in the notation of that paper,
the singular symmetric $2$-tensors $\frac{dw_i}{x}\frac{dw_j}{x}$ are a basis of {\it smooth} 
sections of a bundle which we denote $S^2_0(M)$. Thus, for example, any conformally compact
with smooth conformal compactification is an element of $\calC^\infty S^2_0(M)$.
This normalization differs from the ones in \cite{GL} and \cite{L} and shifts by $2$
the indicial roots. This explains the discrepancy with the numerology in those papers. We have
\begin{proposition} 
Assume that $g$ is a conformally compact, so that $g= x^{-2}\, \olg$, and assume also that 
$|dx|_{\olg} =1$ on $\del M$. Then the set of  indicial roots of ${\cal L}_g$ consists of the pairs:
\[
\zeta_1^\pm = \frac12  \, (n \pm \sqrt{n^2 + 8n}),  \qquad \qquad  \zeta_2^\pm = \frac12  \, (n \pm \sqrt{n^2 + 4n+4}), 
\]
\[
\mbox{and}\qquad \zeta_3^\pm = \frac12 \, (n \pm n) = 0,n. 
\]
Setting $\mu_- = \sup_j \zeta_j^-$ and $\mu_+ = \inf_j \zeta_j^+$, then we have the important inequality
\[
\mu_- : = 0 \ < \ \mu_+ : = n.
\]
\label{pr:3}
\end{proposition}
As already mentioned, the computation of the indicial roots can be found in Lemma 7.1 and 
Lemma 7.5 of \cite{L}. However, to be explicit, the action of 
\[
L_{g_0} : = (\nabla^{g_0})^* \nabla^{g_0} - 2 \,  \stackrel{\circ}{R^{g_0}},
\]
on trace-free symmetric $2$-tensors
\[
h: =  h_{00} \, \frac{ds^2}{s^2} + h_{0i}\, \frac{ds}{s} \, \frac{du_i}{s} + h_{ij} \, 
\frac{du_i}{s}\, \frac{du_j}{s},
\]
is given by
\[
\begin{array}{rllll}
L_{g_0} \, h & = & \displaystyle \left( - \Delta_{g_0} \ h_{00} + 2 n  \ h_{00} - 4 s \ \del_{u_i} h_{0i}\right) \, \frac{ds^2}{s^2} \\[3mm] 
&+ & \displaystyle \left( - \Delta_{g_0} \ h_{0i} + (n+1) \ h_{0i} - 2  s \ \del_{u_j} h_{ij}\right) \, \frac{ds}{s} \, \frac{du_i}{s}  \\[3mm]
& + & \displaystyle \left( - \Delta_{g_0} \ h_{ij} + 2 s \ (\del_{u_j}  h_{0i} + \del_{u_i} h_{0j}) - 2 \ h_{00} \ \delta_{ij} \right) \,  \frac{du_i}{s} \, \frac{du_j}{s} , 
\end{array}
\]
where $\Delta_{g_0}$ is the Laplace-Beltrami operator on hyperbolic space defined in  (\ref{eq:lhs}). 
It is straightforward to check that the indicial roots corresponding to the normal part (i.e.\ $h_{00}$) 
are $\zeta_1^\pm$, the indicial roots corresponding to the mixed part (the $h_{0i}$ terms) are 
$\zeta_{2}^\pm$, and the indicial roots corresponding to the tangential part (the $h_{ij}$ terms) 
are $\zeta_{3}^\pm$. Finally, the action of $L_{g_0}$ on pure-trace symmetric $2$-tensors is given by
\[
L_{g_0} (\tau \, g_0) = (- \Delta_{g_0} \tau + 2n \ \tau) \, g_0,
\]
and so the indicial roots here are also given by $\zeta_{1}^\pm$.

\subsubsection{Function spaces}
Let us now recall the scale of weighted scale invariant H\"older spaces $x^\mu \Lambda^{k,\alpha}_0(M)$.
(These do not provide the optimal boundary regularity for this problem, but they are sufficient for 
our goals here.) For simplicity, definitions will be stated primarily for functions,
but they transfer easily to sections of vector bundles, as indicated briefly below.
We refer to \cite{Ma1} for further discussion and proofs, cf.\  also \cite{L} and \cite{A3}. 

We first define $\Lambda^{0,\alpha}_0(M)$ as the natural `geometric' H\"older space associated
to any fixed smooth conformally compact metric $g = x^{-2} \, \olg$.  Thus if $w=(x,y)$ is a smooth 
coordinate chart near $\del M$, then this space is the closure of $\calC^\infty(\Mbar)$ 
with respect to the norm 
\begin{eqnarray*}
||u||_{0,\alpha} & := & \sup_B \sup_{w,w' \in B} 
\frac{|u(x,y) - u(x', y')|(x + x')^\alpha}{|x -x'|^\alpha + |y - y'|^\alpha} \\
& \cong & \sup_B \sup_{w,w' \in B} \frac{|u(x,y) - u(x', y')|}{\mbox{dist}_{g'} 
\,(w, w')^\alpha} \quad < \ \infty\, ;
\end{eqnarray*}
this supremum is taken first over all points $w=(x,y)$, $w'=(x', y')$,  $ w \neq w'$, lying in some 
coordinate cube $B$ centered at a  point $w_0 = (x_0,y_0)$ of sidelength $\frac12 x_0$, and then over 
all such cubes. We could equivalently replace these cubes by geodesic balls (with respect to $g$)
of radius $1$. Continuing on, we define $\Lambda^{k,\alpha}_0(M)$ to consist of all functions $u$
such that $(x\del_x)^j(x\del_y)^\lambda u \in \Lambda^{0,\alpha}_0(M)$ for all $j+|\lambda| \leq k$.
Noting that if we use the vector fields and $1$-forms $x\del_{w_i}$, $dw_j/x$ and their
tensor products as generators for the sections of all of the tensor bundles, then $\nabla^{g}$ 
involves only differentiations with respect to  $x \, \del_{x}$, $x \, \del_{y_j}$. Hence 
the definitions of these function spaces extend naturally to sections of any of these bundles. 

These norms respect the natural scale invariance of uniformly degenerate operators: in fact, 
for functions $u$ supported in a coordinate chart near the boundary, the norms of $u(w)$ and 
$u_{\e}(w) = u(w/\e)$  are the same. 

For $\mu \in \RR$, define also 
\[
x^\mu\Lambda^{k,\alpha}_0(M) := \left\{u = x^\mu \tilde{u}\ : \quad   
\tilde{u} \in \Lambda^{k,\alpha}_0(M)\right\},
\]
with the corresponding norm denoted $||\cdot||_{k,\alpha,\mu}$.

We shall also have occasion to use weighted $L^2$ spaces. Using any standard 
coordinate chart near the boundary, we define the main `reference' $L^2$ space
\[
L^2_b(M) = L^2\left(M; \frac{dx\,dy}{x}\right)
\]
and its weighted versions
\[
x^\delta L^2_b(M) = \{v =  x^\delta \tilde{v} \ : \quad \tilde{v} \in L^2_b(M)\}.
\]
Notice that the most natural one of these is 
\[
L^2(M;dV_g) = L^2\left(M;\frac{dx\,dy}{x^{n+1}}\right) = x^{n/2}L^2_b(M).
\]
We use $L^2_b(M)$ as the reference space because of the simple numerology that 
\[
x^\mu \Lambda^{k, \alpha}_0 (M) \hookrightarrow x^\delta L^2_b (M)
\]
if and only if $\mu > \delta$. 

\begin{definition} 
The closure in the norm $||\cdot||_{k,\alpha,\mu}$  of the space of smooth symmetric $2$-tensors $\calC^\infty S^2_0(M)$ is denoted $x^\mu \Lambda^{k,\al}_0 \, S^2_0(M)$. 
\end{definition}
In particular, $k \in x^\mu \Lambda^{k,\al}_0 \, S^2_0(M)$ can be decomposed as 
\[
k =  \sum_{i,j=0}^n k_{ij} \, \frac{dw_i}{x} \, \frac{dw_j}{x},
\]
where $k_{ij}\in x^\mu \Lambda^{k,\al}_0 \,(M)$. Similarly, we can define the spaces  $L^2_b \, S^2_0 (M)$ 
and  $L^2 \, S^2_0 (M; dV_g)$ of sections of symmetric $2$-tensors.

\subsubsection{Mapping properties}
Assume that $g$ is a conformally compact metric, so that it can be  written $g = x^{-2}\, \olg$. 
Further assume that  $|dx|_{\olg} =1$. It follows immediately from the definitions that 
\begin{equation}
L_g: x^\mu \, \Lambda^{\ell+2,\alpha}_0 \,S^2_0(M) \longrightarrow x^\mu \, 
\Lambda_0^{\ell,\alpha} \, S^2_0(M)
\label{eq:map}
\end{equation}
is bounded for any $\mu \in \RR$ and $\ell \geq 0$. However, this map does not have closed range when $\mu$ is equal to one of the indicial roots $\zeta_j^\pm$ of $L_g$. This stems from the fact that when $\mu$  is an indicial root, the equation $I(L_g)u = s^\mu$ has solution  $u = c \, s^{\mu}(\log s)$, where $c$ is a constant tensor, and this misses lying  in $s^\mu \Lambda^{2,\alpha}_0$ on account of the logarithmic factor, although $s^\mu$ does lie in this space. 

\begin{proposition}[\cite{Ma1}, cf.\ also \cite{L}] 
If $\mu \in (0,n)$, then  the mapping (\ref{eq:map}) is Fredholm of index zero.
\end{proposition}  
When $\mu$ is not in this range, and is not equal to an indicial root, then  (\ref{eq:map}) still has closed range and is semi-Fredholm, but has either  kernel or cokernel which is infinite dimensional. The operator $L_g$ has analogous  mapping properties when acting between weighted Sobolev spaces, and these are  actually more elementary than the one here between weighted H\"older spaces,  cf.\ \cite{Ma1}. 

The fact that the Fredholm index of $L_g$ is zero when $\mu \in (0,n)$ is proved  using the self-adjointness of $L_g$ on $L^2 \, S^2_0  (M;dV_g)$. For the local deformation theory for Poincar\'e-Einstein  metrics, it is crucial to know whether (\ref{eq:map}) is surjective at  a weight $\mu$ in this interval. This proposition shows that in this  Fredholm range, surjectivity is equivalent to injectivity at the same weight.  This motivates why it is necessary to study the nullspace of $L_g$ more closely. We first state a basic regularity theorem:
\begin{proposition}
Assume that $g$ is a Poincar\'e-Einstein metric. Suppose that  $L_g k = 0$ and that $k \in x^\mu \Lambda^{2,\alpha}_0 \, S^2_0 (M)$ or  $k \in x^\mu L^2_b \, S^2_0(M)$ for some $\mu > 0 $. Then $k$  is polyhomogeneous, and has an expansion 
\[
k \sim \sum k_{j,\ell,r}(y) \, x^{\zeta_j^+ + \ell}(\log x)^r,  \qquad k_{j,\ell,r}  \in \calC^\infty \, S^2_0(\del M),
\]
where the $\zeta_j^+$ are the indicial roots of ${\cal L}_g$ for which $x^{\zeta_j^+} \in  x^\mu \Lambda^{2,\al}_0 (M)$ or  $x^{\zeta_j^+} \in  x^\mu \, L^2_b \, (M)$).  Consequently,  for any such $k$ we necessarily have $|k|_{g} = \calO(x^{n})$. 
\label{pr:5}
\end{proposition}
Note in particular that this result applies to any $k$ which belongs  to $\ker L_g \cap L^2 \, S^2_0(M;dV_g)$. We remark also that elements of  $\ker L_g$ which lie in a weighted space with $\mu \leq 0$ are no longer necessarily polyhomogeneous, and their precise regularity is  determined by the regularity of their leading asymptotic coefficients,  c.f.\ \cite{Ma1}. 

An immediate consequence of this result is that the nullspace of $L_g$  on any one of the $x^\mu 
L^2_b \, S^2_0 (M)$ or $x^\mu \Lambda^{2,\al}_0 \, S^2_0 (M)$ for $\mu \in (0,n)$ are all the same. 
\begin{definition} 
A Poincar\'e-Einstein metric $g$ is said to be nondegenerate if the nullspace of $L_g$ on 
$L^2 \, S^2_0(M;dV_g)$ is trivial.
\label{de:ndg}
\end{definition}

The local deformation theory for \PE metrics is simplest when the $L^2$  nullspace for $L_g$ is 
trivial. The triviality of this nullspace is verified, in increasingly general settings, by 
Graham and Lee \cite{GL}, Biquard \cite{Bi} and Lee \cite{L}. Anderson modifies this  
approach by including the variation of the conformal infinity as an explicit variable, and shows 
that the resulting map is always surjective. We explain this more carefully now.

Choose a smooth boundary defining function $x$ according to Lemma~\ref{le:1}. The flow lines 
for $\nabla^{\olg} x$ determine a product structure $[0,x_0)_x \times \del M$ near the boundary,
and we let $\pi$ be the projection onto the second factor. Fix a smooth cutoff function 
$\chi(x)\geq 0$ which equals $1$ for $x \leq x_0/2$ and which vanishes for $x \geq x_0$, 
and define the extension map 
\[
e: \calC^{2,\al}S^2(\del M)  \longrightarrow  x^{-2}\calC^{2,\al}S^2(M) \subset \Lambda^{2,\al}_0 S^2_0(M)
\]
\[
e(\eta) \quad =  \quad \frac{\chi(x) \, \pi^*(\eta)}{x^2}.
\]
When $\eta$ is small, 
\[
g_\eta : =  g + e(\eta)
\]
is a conformally compact metric. Proposition 1.4.6 in \cite{Bi} asserts that for any 
$\eta \in {\calC}^{2,\al}\, S^2(\del M)$ with small norm, there exists a unique 
diffeomorphism $\varphi_\eta$ close to the identity such that the metric
\begin{equation}
\tilde g_\eta : =  {\varphi_\eta}^* \, g_\eta
\label{eq:tg}
\end{equation}
satisfies the Bianchi gauge condition 
\begin{equation}
B^{g}(\tilde g_\eta) = 0.
\label{eq:jh}
\end{equation}
In particular, $\varphi_0 =Id$. In fact, the construction yields that
$\varphi_\eta \in {\mathfrak D}^{3,\al}_2(M)$, the group of $\calC^{3,\al}$ diffeomorphisms 
induced by vector fields in $x^{2}\,\Lambda^{3,\al}_0(M)$. Note that $\varphi_\eta$ restricts 
to the identity map on $\del M$. 

The mapping 
\[
\begin{array}{rlllll}
\Gamma = \Gamma^g : {\calC}^{2,\al}\, S^2(\del M) & \longrightarrow & \Lambda^{2,\al}_0 \, S^2_0(M)\\[3mm]
\eta & \longmapsto & \tilde g_\eta 
\end{array}
\]
is smooth and has differential
\[
\left. D\Gamma \,  \right|_{\eta=0}(\gamma) = e(\gamma) + h(\gamma),
\]
where 
\[
h:{\calC}^{2,\al}\,S^2(\del M) \longrightarrow  x^2 \Lambda^{2,\al}_0 S^2_0 (M)
\]
is bounded. Differentiating (\ref{eq:jh}) with respect to $\eta$ yields
\begin{equation}
B^{g}(e(\gamma) + h(\gamma)) = 0.
\label{eq:jhh}
\end{equation}

Following Anderson \cite{A3}, define the nonlinear mapping
\begin{equation}
\tilde N_g(\eta, k) := \Ric^{\tilde g_\eta + k} + n(\tilde g_\eta + k ) + 
(\delta^{\tilde g_\eta+k})^* B^{\tilde g_\eta} (k).
\label{eq:mainext}
\end{equation}
This maps a neighbourhood of the origin in ${\calC}^{2,\al}\,S^2(\del M) \, 
\oplus \, x^\mu \Lambda^{2,\al}_0 \, S^2_0(M)$ into $x^\mu \Lambda^{0,\al}_0 \, S^2_0 (M)$, 
provided $0 < \mu \leq 2$.  Indeed, it follows from  (\ref{eq:esterr}) together with the choice of 
defining function $x$ as in Lemma~\ref{le:2} that $|\mbox{Ric}^{{\tilde g}_\gamma} + n \, 
{\tilde g}_\gamma|_g  \in x^2 \Lambda^{0,\al}_0 S^2(M)$. 

Now
\[
\begin{array}{rllllll}
2 \, D {\tilde N}_g|_{(0,0)}(\gamma, \kappa)  & := & {\tilde L}_g (\gamma, \kappa) \\[3mm]
&  = & L_g (e(\gamma)  + h(\gamma) + \kappa) - 2 \, (\delta^{g})^*  \, B^g (e(\gamma) + h(\gamma)).
\end{array}
\]
Using (\ref{eq:jhh}) to simplify this, we obtain finally
\begin{equation}
\tilde L_g (\gamma, \kappa) = L_g (e(\gamma)  + h (\gamma) + \kappa).
\label{eq:equ}
\end{equation}
\begin{proposition}[\cite{A3} Theorem 1.2] 
Fix a \PE metric $g$ on $M$ and let $m = \dim \ker L_g$. Then there exists an $m$-dimensional 
subspace $\calS \subset \calC^{2,\al}\,S^2_0(\del M)$, such that 
\begin{equation}
\tilde L_{g}: \calS \oplus x^\mu \, \Lambda^{2,\al}_0\, S^2_0(M) \longrightarrow x^\mu \, 
\Lambda^{0,\al}_0 \, S^2_0(M)
\label{eq:exmapfd}
\end{equation}
is surjective for any $0 < \mu \leq 2$. 
\label{pr:qqq}
\end{proposition}
The proof boils down to showing that for every $\kappa \in \ker L_g$, there exists a 
$\gamma \in \calC^{2,\al}(\del M)$ such that 
\begin{equation}
\int_M \langle \tilde L_g \, (\gamma ,0), \kappa \rangle \, dV_g 
= \int_M \langle L_g \, (e(\gamma)+h(\gamma) ), \kappa \rangle \, dV_g \neq 0.
\label{eq:dfg}
\end{equation}
Integrating by parts, this is the same as the requirement that
\begin{equation}
\lim_{a \rightarrow 0} \int_{\del M_a} \left( \langle \nabla^g_N \gamma, 
\kappa \rangle - \langle \gamma ,  \nabla^g_N \kappa \rangle\right) \, dV_g \neq 0;
\label{eq:nnn}
\end{equation}
here $M_a = \{x \geq a\}$ and $N = \nabla^g x/ |\nabla^g x|$ is the unit normal to $\del M_a$. 
Using that $\kappa$ is polyhomogeneous with leading term $x^n$ and expanding the left side of 
(\ref{eq:nnn}) as $a \searrow 0$, we deduce that (\ref{eq:nnn}) fails to hold precisely
when $\gamma$ lies in a codimension $1$ subspace $V(\kappa)$. The intersection of these
subspaces as $\kappa$ varies over $\ker L_g$ has codimension $m$, which finishes the proof.

As a corollary, the implicit function theorem now implies that the moduli space $\calE(M)$ 
of \PE metrics on $M$ is always a smooth Banach manifold.  It also follows that the restriction 
of the projection $\Pi$ to  $\calE$ is Fredholm of index zero; by the Sard-Smale theorem, its 
image is (at worst) a (Banach) variety of finite codimension (provided  $\calE(M) \neq \emptyset$).  

In our gluing construction it will be convenient to work with variations $\gamma$ and
diffeomorphisms $\varphi$ which are supported away from the points $p_j \in \del M_j$. We now
argue that Proposition \ref{pr:qqq} and its proof are stable under small perturbations. 

Fix a boundary coordinate chart $(x,y)$ centered at $p \in \del M$. Letting $\psi(r)$ be a smooth
cutoff function which vanishes for $r \leq 1/2$ and equals $1$ for $r \geq 1$, and $\psi_\tau(r): = 
\psi(r/\tau)$, we define
\[
\hat \varphi_\eta : = (1 -\psi_\tau(r)) \, \mbox{Id} + \psi_\tau(r) \, \varphi_\eta,
\qquad r = |(x,y)|, 
\]
in that chart and $\hat\varphi_\eta = \varphi_\eta$ elsewhere; we also define
\[
\hat g_\eta : =    \hat {\varphi_\eta}^* \, (g + \psi_\tau e(\eta)).
\]
Clearly $\hat g_\gamma = g$ on the ball of radius $\tau$ around $p$. 

Replacing $\tilde g_\eta$ by $\hat g_\eta$ in the definition of $\tilde N_g$ yields 
a new nonlinear operator $\hat N_g$, with linearization $\hat L_g$. 
We use $\varphi_\gamma \in {\mathfrak D}^{3,\alpha}_2(M)$ to show that 
$(\tilde L_g - \hat L_g) (\gamma, 0) \to 0$ as $\tau \to 0$. Now suppose that (\ref{eq:nnn}) 
holds for some $\kappa$ and $\gamma$. To prove that
\[
\lim_{a\rightarrow 0} \int_{\del M_a} \left( \langle \nabla^g_N (\psi_\tau \gamma),  \kappa 
\rangle - \langle \psi_\tau \, \gamma ,  \nabla^g_N \kappa \rangle\right) \,  dV_g \neq 0
\]
when $\tau$ is small, we only need observe that the commutator $[\nabla^g_N, \psi_\tau]$
is bounded and supported on an increasingly small set, hence its contributionto the 
integral becomes increasingly negligible. This establishes the analogue of Proposition~\ref{pr:qqq}. 

\section{The approximate solution}

We now commence with the construction.  Suppose that $(M_j,g_j)$,  $j=1,2$,  are Poincar\'e-Einstein, $\dim M_j = n+1$. Fix points  $p_j \in \del M_j$ and, as in the introduction, consider the boundary  connected sum $M_1 \#_b M_2$, obtained by excising half-balls around  the points $p_j$ and identifying their hemispherical boundaries. We  describe in the next paragraphs how to perform this operation on the  scale of a small parameter $\e$, and shall denote the resulting  manifold $M_\e$. 

We fix a defining function $x$ as in Lemma~\ref{le:2} so that 
\[
g_1 = \frac{dx^2 + h^{(1)}(x)}{x^2}
\]
where $h^{(1)}(x)$ is a $\calC^{2,\al}$ family of $\calC^{2,\al}$ metrics on  $\del M_1$ and
\begin{equation}
h^{(1)} = h^{(1)}_0 + {\calO} (x^2),
\label{eq:est1}
\end{equation}
where $h^{(1)}_0 : = h^{(1)}(0)$ is a representative of $\frakc (g)$  the conformal infinity of $g_1$. We further fix Riemann normal  coordinates $y$ centered at $p_1 \in \del M$ for the metric  $h^{(1)}_0$. By definition of normal coordinates, 
\begin{equation}
h^{(1)}_0 = dy^2 + {\mathcal O}(|y|^2).
\label{eq:est2}
\end{equation}
Then, in the boundary normal coordinates $w=(x,y)$ we have
\begin{equation}
g_1  =  \frac{dx^2 + dy^2}{x^2} +  \, \sum_{\alpha,\beta = 1}^n k_{\alpha \beta}^{(1)} \, \frac{dy_\alpha}{x} \,  \frac{dy_\beta}{x} .
\label{eq:g1nc}
\end{equation}
Thus 
\[
k^{(1)} : =  \sum_{\alpha,\beta = 1}^n k_{\alpha \beta}^{(1)}(w) \, 
\frac{dy_\alpha}{x} \, \frac{dy_\beta}{x},
\]
is a $2$-tensor which measures the discrepancy of $g_1$  from the standard hyperbolic metric $g_0$ and (\ref{eq:est1}) together with (\ref{eq:est2}) give the estimate for the coefficients of the discrepancy tensor
\begin{equation}
k_{\alpha \beta}^{(1)} (w) = \calO(|w|^2).
\label{eq:measdisc}
\end{equation}
Similarly, in terms of boundary normal coordinates $w'=(x',y')$ near  $p_2$  on $M_2$, we can decompose $g_2 = g_0 + k^{(2)}$, and the  coefficients of this discrepancy tensor (relative to the coframe  $dw_j'/x'$) are $\calO(|w'|^2)$. 
 
Let $A_\e$ and $A_\e'$ denote the annuli $\{\e/2 \leq |w| \leq 2\e\}$ and $\{\e/2 \leq |w'| \leq 2\e\}$, respectively. Identifying these by  means of the inversion mapping $w' = I_\e(w)$ where   $I_\e(w) : = \e^2 w/|w|^2$, we define the smooth manifold with  boundary 
\[
M_\e = \bigg(M_1 - B_{\e/2}(p_1)\bigg) \, \bigcup_{I_\e} \,  \bigg(M_2 - B_{\e/2}(p_2)\bigg).
\]
Note that the annulus $A_\e \sim A_\e'$ is naturally embedded in  $M_\e$.

It will be more convenient to use a rescaling of these coordinate  systems, so that we may regard the gluing region as a fixed annulus  $A$. Thus define the dilation $R_\e$, which sends $w$ to $\e w$ (and  $w'$ to $\e w'$). The annuli $A$ and $A'$ of inner and outer radii  $1/2$ and $2$ in the $w$, $w'$ coordinates are mapped by $R_\e$ to  $A_\e$ and $A_\e'$, respectively, and are identified by the fixed  inversion $I(w) = w/|w|^2$. Observe that  $I = R_{\e}^{-1} \, I_\e \, R_\e$.

The metrics $g_{j,\e} = R_\e^*(g_j)$ are defined on the half-ball of radius $C/\e$ for some $C>0$; these are just isometric forms of the  initial metrics $g_j$. We define a family of conformally compact  metrics $g_\e$ on $M_\e$ by pasting these together. Thus let $\chi(r)$ be a nonnegative, smooth cutoff function which  equals $1$ for $r = |w| \geq 2$ and vanishes for $r \leq 1/2$. Then set
\[
g_\e = \chi(r) g_{1,\e} + (1-\chi(r))I^*(g_{2,\e}).
\]
This is our approximate solution; it is a conformally compact metric  on $M_\e$ and agrees with the original metrics $g_1$ and $g_2$ outside  of the half-balls $B_{2\e}(p_j)$.

We now estimate the discrepancy of $g_\e$ from being Einstein. To this end, observe that 
\[
g_{1,\e} = \frac{dx^2 + dy^2}{x^2} + \sum_{\alpha,\beta = 1}^n  \left(R_\e^* k_{\alpha \beta}^{(1)}
\right)(w) \frac{dy_\alpha}{x}\frac{dy_\beta}{x},
\]
since $R_\e^*(dw_j/x) = dw_j/x$. Clearly
\begin{equation}
|R^*_\e k^{(1)}_{\alpha \beta}(w)| \leq C \e^2, 
\label{eq:errk1}
\end{equation}
uniformly for $w \in A$, with $C$ independent of $\e$. The  coefficients $R^*_\e k^{(2)}$ satisfy the same estimate for $w' \in A'$.  To compute $I^*(g_{2,\e})$, it suffices to concentrate on the term  $I^* \, R^*_\e \, k^{(2)}$  since $I$ is an isometry of $g_0$. We have
\begin{equation}
I^{\, *} \left(\frac{dy_\alpha'}{x'}\right) = \frac{dy_\alpha}{x} - 2 \frac{y_\alpha}{r}\left( \frac{x}{r}\frac{dx}{x} +  \sum_{\beta=1}^n \frac{y_\beta}{r}\frac{dy_\beta}{x} \right),
\label{eq:1form}
\end{equation}
where $r = (x^2 + |y|^2)^{1/2}$. This then gives for all $w\in A$
\begin{equation}
\begin{array}{rllll}
I^{\,*} \, R^*_\e \, k^{(2)} & = & \displaystyle \calO(\e^2 \, x^2)\frac{dx^2}{x^2}  +   \displaystyle \sum_{\alpha=1}^n \calO(\e^2 \, x) \frac{dx}{x} \frac{dy_\alpha}{x} \\[3mm]
&+& \displaystyle \tilde{I}^{\, *}  R_\e^*  \, k^{(2)} +  \sum_{\alpha, \beta=1}^n \calO(\e^2 \, x^2)  \frac{dy_\alpha}{x}\frac{dy_\beta}{x} ,
\end{array}
\label{eq:g2err}
\end{equation}
where $\tilde{I}(y) = y/|y|^2$ is the restriction of the inversion $I$  to the boundary. 
The expression for the final term uses that for $w \in A$ and $x \to 0$, $r = |y|+ {\calO}(x^2)$.

Note that the first two terms and the last term on the right in (\ref{eq:g2err}) vanish  at $\del M_\e$.  Hence the conformal infinity of $g_\e$ is represented  by the metric $h_{0,\e}$ which is obtained by identifying the annuli   $1/2 \leq |y| \leq 2$, $1/2 \leq |y'| \leq 2$ in the rescaled normal coordinates on $\del M_1$ and $\del M_2$ using the  inversion $\tilde{I}$, and pasting together the metrics $h_0^{(1)}$  and $h_0^{(2)}$ with the cutoff function $\chi(|y|)$. This will be  important in \S 6 when we discuss the Yamabe type of $\frakc(g_\e)$. 

The expansion (\ref{eq:g2err}) has many other consequences. Observe  that, in the annulus $A$, the metric $g_\e$ can be expanded as
\begin{equation}
g_\e =  \displaystyle \left( 1 + {\calO} (\e^2 x^2) \right) \,  \frac{dx^2}{x^2} + \sum_{\alpha=1}^n {\calO} (\e^2 x) \, \frac{dx}{x}   \frac{dy_\alpha}{x} + \displaystyle \sum_{\alpha , \beta=1}^n \left(  \delta_{\alpha\beta} + {\calO} (\e^2) \right) \,\frac{dy_\alpha}{x}   \frac{dy_\beta}{x}
\label{eq:expg}
\end{equation}
In particular, this implies that
\begin{equation}
|dx|^2_{x^2 \, g_\e} = 1+ {\calO} (\e^2 \, x^2),
\label{eq:dx}
\end{equation}
in $A$.

It remains to estimate $N_{g_\e}(0)$, which measures the discrepancy 
of $g_\e$ from being Einstein (and in the proper gauge). By definition of (\ref{eq:main}) we have
\[
N_{g_\e}(0) = \Ric^{g_\e} + ng_\e.
\]
This is supported in $A$, and since both $R_\e^* k^{(1)}$ and
$I^* R_\e^* k^{(2)}$ are $\calO(\e^2)$, along with their derivatives,
this error term is also $\calO(\e^2)$ in $A$. There is an improved
estimate as $x\to 0$. Indeed, if $\olg_\e := x^2 g_\e$, then $\Ric^{\olg_\e} = \calO(\e^2)$ in $A$. 
Furthermore, if $f= - \log x$, then
\[
\nabla^{\olg_\e} df  - df \circ df = \frac{1}{x} \, {\calO} (\e^2), \qquad \mbox{and}
\qquad \Delta_{\olg_\e} f - |df|_{\olg_\e} = \frac{1}{x} \, {\calO} (\e^2)
\]
in $A$. Now use (\ref{eq:ricg}) with $\olg_\e : = x^2 \, g_\e$  and $f : = -\log x$ to conclude that  
\[
\Ric^{g_\e} + n g_\e = \frac{1}{x} \, {\calO} (\e^2),
\]
in $A$. In particular, taking the norm with respect to $g_\e$, we obtain finally the
\begin{proposition}
For the metric $g_\e$ on $M_\e$, the tensor  $N_{g_\e}(0) = \Ric^{g_\e} + n g_\e$ vanishes outside the annulus $A$;  in $A$ its pointwise norm with respect to $g_\e$ satisfies 
\[
|N_{g_\e}(0)|_{g_\e} \leq C \, \e^2 \, x,
\]
where $C$ is independent of $\e$. 
\label{pr:6}
\end{proposition}

\section{Linear estimates}

Let $L_{g_\e}$ denote the linearization of the map $k \to N_{g_\e}(k)$ at $k=0$. Then
\[
L_{g_\e} = L_{0\e} + B_\e,
\]
where
\[
L_{\e} : = (\nabla^{g_\e})^* \nabla^{g_\e} -  2 \, \stackrel{\circ}{R^{g_e}}, 
\]
and
\[
B_\e (\kappa) : = \Ric^{g_\e}\circ \kappa + \kappa \circ \Ric^{g_\e} +   2 \, n \, \kappa .
\]
Of course $B_\e$ would vanish if $g_\e$ were Einstein, and in any  case, $B_\e$ is supported in $A$ and has coefficients which are  $\calO(\e^2)$.  Our goal in this section is to verify a certain  weighted estimate for $L_{g_\e}$, which we now explain.

In the rescaled coordinates $w$ and $w'$, $M_\e$ contains an expanding  annular region $\calT_\e = \{C\e \leq |w| \leq C/\e\}$; the outer  portion $1/2 \leq |w| \leq C/\e$ corresponds to a region in $M_1$ and  the inner portion $2 \geq |w| \geq C\e$ corresponds by inversion to a  region in $M_2$. We introduce polar coordinates 
\[
(r,\phi,\omega) \in {\mathbb R}^+ \times [0, \pi/2] \times S^{n-1},  \qquad  w := (x,y)= (r\, \cos \phi, r \,\sin \phi\, \omega) .
\]
Setting $s = \log r$, and dropping an irrelevant additive constant, then 
\[
\calT_\e = \{(s,\phi,\omega): - s_\e \leq s \leq s_\e\}, \qquad \mbox{where}\quad s_\e = -\log \e.
\]
As $\e \to 0$, $\calT_\e$ expands to fill out $\calT_0 = \RR \times S^n_+$, and the metric $g_\e$ converges (in $\calC^\infty$ on compact sets) to 
\[
\frac{1}{\cos^2\phi} \, \left(ds^2 + d\phi^2 + \sin^2 \phi \, d\omega^2 \right),
\]
This is the standard hyperbolic metric $g_0$ on $\HH^{n+1}$, written in warped product form; the hemisphere $S^n_+$ with metric $(d\phi^2 + \sin^2  \phi  \, d\omega^2)/\cos^2 \phi$ is isometric to $\HH^n$. Notice that $L_{g_\e}$ and $L_{\e}$ both converge in this central region to $L_{g_0}$, the linearized Einstein operator on $\HH^{n+1}$. 
n
Define a weight function $w_\e$ on $M_\e$ which is a smoothed version 
of the function
\[
\tilde{w}_\e(s) = \left\{
\begin{array}{rcl}
& \cosh s/\cosh s_\e \quad & \mbox{in}\  \calT_\e, \\[3mm]
& 1 & \mbox{in}\ M_\e - \calT_\e.
\end{array}
\right.
\]
(We require that $w_\e$ is smooth, and agrees with $\tilde{w_\e}$ except on a small neighborhood of $\del \calT_\e$, where it remains bounded between $1/2$ and $2$.) Now let $\rho$ denote a boundary  defining function for $M_\e$ which is $\calC^\infty$ and agrees (up to  a small smoothing near $\del B_1(p_j)$) with fixed boundary defining  functions  $\rho_j$ on $M_j - B_{1}(p_j)$ (in {\it unrescaled}  coordinates!)  and with $t = \cos \phi$ in $\calT_\e$. In terms of  these, we define the doubly weighted H\"older spaces
\[
\rho^\mu (w_\e)^\nu \Lambda^{k,\al}_0(M_\e),
\]
which contain all functions of the form  $u = \rho^\mu \, (w_\e)^\nu \, v$, with  $v \in \Lambda^{k,\al}_0 (M_\e)$. We define the space of symmetric  $2$-tensors $\rho^\mu (w_\e)^\nu \Lambda^{k,\al}_0S^2_0(M_\e)$  similarly. We denote the corresponding norm  $||\cdot ||_{k,\al,\mu,\nu}$.  

\begin{proposition}
Suppose that neither of the operators $L_{g_1}$ nor $L_{g_2}$ has a  nontrivial $L^2$ nullspace. Also, let $0 < \mu = \nu < n/2$. Then for  $\e$ sufficiently small, the operator
\[
L_{g_\e} : \rho^\mu (w_\e)^\mu \Lambda^{2,\al}_0 S^2_0(M_\e)  \longrightarrow \rho^\mu (w_\e)^\mu \Lambda^{0,\al}_0 S^2_0(M_\e) 
\]
is an isomorphism and its inverse $G_{g_\e}$ has norm bounded  independently of $\e$. 
\label{pr:7}
\end{proposition}
{\bf Proof:} Suppose this proposition were false. Then there would  exist a sequence $\e_j \to 0$  and sequences of $2$-tensors $h_j$ for  which
\[
||h_j||_{2,\al,\mu,\mu} = 1,
\]
while 
\[
||L_{g_{\e_j}} h_j||_{0,\al,\mu,\mu} \to 0 .
\]
Rewriting these norm inequalities gives the estimates
\begin{equation}
|h_j(z)| \leq (\rho(z) \, w_{\e_j}(z))^\mu, \qquad \mbox{and}\qquad |L_{g_{\e_j}} h_j(z)| \leq \eta_j (\rho(z) \, w_{\e_j}(z))^\mu,
\label{eq:bounds}
\end{equation}
for all $z \in M_{\e_j}$, where $\eta_j \to 0$. We shall use a blow-up analysis to show that this leads to a contradiction.

Suppose that the supremum of the pointwise norm  $\rho^{-\mu}(w_\e)^{-\mu} |h_j|$ occurs at some point $q_j$. (If this  supremum is not attained anywhere in the interior of $M_\e$, then it  is enough to assume that the value of this function at $q_j$ is larger than  half its supremum.) Possibly passing to a subsequence, there are  several cases which may arise:
\begin{enumerate}
\item[(i)] $q_j$ converges to a point $q$ in the interior of $M_1$ or  $M_2$;
\item[(ii)] $q_j$ lies in $\calT_{\e_j}$ for every $j$ and its coordinates  $(s_j,t_j,\omega_j)$ (where $t_j = \cos \phi_j$) satisfy  $|s_j| \leq C$, $t_j \geq c > 0$;
\item[(iii)] $q_j$ lies in $\calT_{\e_j}$ for every $j$, and its  coordinates $(s_j,t_j,\omega_j)$  satisfy $|s_j| \to \infty$,  $t_j \geq c > 0$;
\item[(iv)] $q_j$ converges to a point $q$ on $\del M_1 - \{p_1\}$  or $M_2 - \{p_2\}$; 
\item[(v)] $q_j$ lies in $\calT_{\e_j}$ for every $j$ and its coordinates  $(s_j,t_j,\omega_j)$ satisfy $|s_j| \leq C$, $t_j \to 0$;
\item[(vi)] $q_j$ lies in $\calT_{\e_j}$ for every $j$ and its coordinates  $(s_j,t_j,\omega_j)$ satisfy $|s_j| \to \infty$, $t_j \to 0$.
\end{enumerate}
These will be ruled out in turn. In each case, we define a new  sequence $k_j$, either by dividing by a normalizing constant so that  $|k_j(q_j)|$ is bounded away from $0$, or else by rescaling the  independent variable, or both. Extracting a subsequence if this is  necessary, we can assume that this new sequence converges to a  symmetric $2$-tensor  $k$ which is defined, either one of the $M_j$ or on $\HH^{n+1}$, and  which is a solution of the linearized Einstein  equation there. The  goal is to show that it satisfies certain $L^\infty$ bounds, and then to show that no such solution exists. In the next several  paragraphs, we deduce the existence of this limiting tensor $k$ and  deduce the bounds it must satisfy. Only afterwards do we show that  these bounds preclude its existence. 

In case (i), $\rho(q_j) \geq c > 0$ and $w_{\e_j}(q_j) = 1$ (or at  least is bounded away from zero), so we take 
\[
k_j := h_j.
\]
Suppose,  for example, that $q_j \to q \in M_1$. Then $k_j$ converges (in  $\calC^\infty$ on compact sets of $M_1$) to a symmetric $2$-tensor $k$  on $M_1$ which satisfies $L_{g_1} \, k = 0$ and also, from  (\ref{eq:bounds}), the bound $|k(z)| \leq \rho_1(z)^\mu$. Furthermore,  $k$ is nontrivial since $k(q) \neq 0$. This case can be ruled out  immediately since we are assuming that the $L^2$ nullspace of  $L_{g_1}$ is trivial, and hence, by Proposition~\ref{pr:5}, so is its nullspace in $\rho_1^\mu \Lambda^{2,\al}_0 \, S^2_0 (M_1)$ for $\mu > 0$. 

Now consider case (ii).  As $\e \to 0$,  $\calT_\e \to \calT_0 = \HH^{n+1}$, and by assumption, $q_j$ remains  in a compact set of $\calT_0$. Thus we may assume that  $q_j\to q \in \calT_0$. We have again that $\rho(q_j) \geq c' > 0$,  whereas $w_{\e_j}(z) \to 0$. However, for $z \in \calT_0$, 
\[
\frac{w_{\e_j}(z)}{w_{\e_j}(q_j)} =\frac{\cosh s}{\cosh s_j} \to c''  \cosh s, \qquad c'' > 0.
\]
(As before, and as in all of the remaining cases, this convergence is  $\calC^\infty$ on compact sets.) Thus if we define 
\[
k_j := (w_{\e_j}(q_j))^{-\mu} \, h_j ,
\]
then $k_j \to k$ where $k(q) \neq 0$, $L_{g_0} \, k = 0$, and  $|k(z)| \leq c\, (\cos \phi \, \cosh s)^\mu$. 

In case (iii), we first recenter the coordinates on $\calT_\e$ by  replacing the independent variable $s$ by $s-s_j$, where $s_j$ is the  $s$ coordinate of $q_j$ (so $|s_j|\to \infty$, by hypothesis), and 
then define 
\[
k_j(s,\phi,\omega) : = (w_{\e_j}(q_j))^{-\mu} \, h_j(s-s_j,\phi,\omega).
\]
But $|h_j(s-s_j,\phi,\omega)| \leq c (w_{\e_j}(s-s_j))^\mu$ and
\[
\frac{w_{\e_j}(s-s_j)}{w_{\e_j}(q_j)} = \frac{\cosh(s-s_j)}{\cosh S_{ \e_j}} \frac{\cosh S_{\e_j}}{\cosh s_j} \to c\, e^s, \qquad c > 0,
\] 
so the limit tensor $k$ is nontrivial, satisfies $L_{g_0} \, k = 0$ and the pointwise bound $|k(z)| \leq c (\cos \phi  \, e^s)^\mu$. 

In case (iv), suppose that $q_j$ remains in some boundary coordinate  chart $(x,y)$ and has coordinates $(x_j,y_j)$, where $x_j \to 0$ and  $y_j \to y_0$.  We may as well assume that $y_0 = 0$, and then define 
\[
k_j(x,y) : = x_j^{-\mu} \, (R_{x_j^{-1}})^* h_j. 
\]
Then $k_j \to k \not\equiv 0$, where $k$ is defined on all of  $\HH^{n+1}$, satisfies $L_{g_0} k = 0$, and the bound  $|k| \leq c \, x^\mu$. Note that this is the exact same bound as in  the previous case, once we have change coordinates  $x = \cos \phi \, e^{s}$.  

In cases (v) and (vi), we can do nearly identical rescalings and  obtain to obtain nontrivial tensors $k$ defined on all of $\HH^{n+1}$  which satisfy $L_{g_0} \, k = 0$ and the bounds  $|k(z)| \leq c\, (\cos \phi \, \cosh s)^\mu$ and   $|k(z)| \leq c (\cos \phi \, e^s)^\mu$, respectively.  

We have already eliminated case (i), so it remains to rule out the  existence of nontrivial solutions of $L_{g_0} \, k = 0$ on $\HH^{n+1}$  which are bounded either by $(\cos \phi \, \cosh s)^\mu$ or  $(\cos \phi \, e^s)^\mu$. Observe that the latter case is included in  the former. It suffices  to check that any such solution lies in  $\rho^\delta L^2_bS^2_0(\HH^{n+1})$ for some $\delta > 0$, where $\rho$ is a global defining function on $\HH^{n+1}$. This is  because by Proposition~\ref{pr:5}, $k$ would then be polyhomogeneous, and in particular decay like $\rho^{n}$ on the entire boundary, and  this would then be ruled out by nondegeneracy. This is a simple  calculation. It suffices to work in the ball $|w| = |(x,y)| < 1$ and,  since $|k(z)| \leq c\, (x/r)^\mu \, r^{-\mu}$, where  $r = (x^2+ |y|^2)^{1/2}$, we have to show that
\[
\int_{x^2+|y|^2\leq 1} \left(\frac{t}{r}\right)^{2\mu} x^{-2\delta}\,  \frac{dx\,dy}{x}  < +\infty
\]
Performing the change of variable $x= r\, \cos \phi$ and  $y=  r\, \sin  \phi\, \omega$, we must show
\[
\int_{t=0}^1 \int_{r=0}^1 t^{2\mu - 2\delta - 1} \, r^{-2\mu-2\delta + n -  1}\, dr \,dt <\infty,
\]
and this holds provided 
\[
\delta < \mu < \displaystyle \frac{n - 2 \delta}{2} .
\]
Since our only other requirements are that $\delta > 0$ and  $0 < \mu < n$, we see that this is easily satisfied if we choose  $\delta$ sufficiently small, provided $\mu \in (0, n/2)$. \hfill  $\Box$

It is fortuitous that all the possible limiting cases which arise here can be handled solely by Proposition~5. However, even if this were  not the case, there is a complete theory for the mapping properties  and regularity of solutions of elliptic uniformly degenerate operators on doubly weighted spaces $x^\mu r^\nu \Lambda^{k,\al}_0(M)$ and $x^\mu \, r^\nu \, L^2_b(M)$, which would have given us the same  kind of conclusions (this would only extend the range in which the  weight parameter $\mu$ can be chosen to $(0,n)$). This more intricate  linear theory is useful in many other problems and we shall return to  it elsewhere.

It remains to adapt this theorem to the case when at least one of the  summands $(M_j,g_j)$ is
degenerate. Recall that, given $\eta^{(j)} \in {\calC}^{2, \alpha} \, S^2_0(\del M_j)$, 
we have defined a metric $\hat g_{\eta^{(j)}}$ on $M_j$. Since these metrics are identically 
equal to $g_j$ near $p_j$, they can be glued together as in \S 3 to produce a metric 
$g_{\e,\eta}$ on $M_\e$,  $\eta = (\eta^{(1)}, \eta^{(2)})$. This allows one to define 
a nonlinear mapping 
\[
\hat N_{g_\e} (\eta, k) : =  \Ric^{\hat g_{\e, \eta} + k} + n(\hat g_{\e, \eta} + k) + 
(\delta^{\hat g_{\e, \eta} +k})^* B^{\hat g_{\e, \eta}}(k),
\]
with linearization at $0$ $\hat L_{g_\e}$.

Let $K_j = \ker L_{g_j}$ and write $m_j : = \dim K_j$. By our current assumption, at least one of the 
$m_j$ is nonzero, and to be definite we suppose that both are. According to Proposition~\ref{pr:qqq}, 
or rather, its modification at the end of \S 2, for $j=1,2$, there exist symmetric $2$-tensors 
$\gamma^{(j)}_i$, $i = 1, \ldots, m_j$, with span $\calS_j$, which are supported away
from $p_j$, and such that the mapping (\ref{eq:exmapfd}) on $M_j$ is surjective. 
We write $\calS = \calS_1 \oplus \calS_2$. 

Recall that $L_{g_\e}$ is the linearization of $\hat N_{g_\e}$ with respect to the 
second factor, i.e., 
\[
L_{g_\e} =  \hat L_{g_\e} (0, \cdot).
\]
This operator is self-adjoint and Fredholm on $L^2S^2_0(M_\e;dV_{g_\e})$, so for some $a > 0$, 
$\mbox{spec}(L_{g_\e}) \cap (-a,a)$ is discrete. The eigenvalues in this range vary continuously 
with $\e$ (so long as they stay within  this interval) but the multiplicity of $0$ as an 
eigenvalue may not be constant. To regain stability we consider the set $\calP_\e$ of all 
`small eigenvalues', which are by definition those which tend to $0$ as $\e \to 0$. We may 
choose $a>0$ so that $\calP_\e$ coincides  precisely with the set of eigenvalues of $L_{g_\e}$ 
lying in $(-a,a)$.  The sum of the corresponding eigenspaces is denoted  $V_\e$ and it 
is standard, cf.\ \cite{Ka}, that $V_\e$ varies continuously with  $\e$. 

Suppose that $h_j \in K_j$, $||h_j||_{L^2} = 1$; using cutoff functions  $\psi(r_j/\tau)$ 
as above, with $r_j = \mbox{dist}\,(w,p_j)$, we define $\tilde{h}_j = \psi(r_j/\tau) h_j$, 
extended by $0$ to the rest of $M_\e$ (i.e. on the other factor $M_{j'}$, $j' \neq j$). 
Then $\tilde{h}_j \to h_j$ in $L^2$ as $\tau \to 0$, and  moreover, if 
\[
\Pi_\e:L^2(M_\e,g_\e) \to V_\e
\]
is the  orthogonal projection, then  $||(I-\Pi_\e)\tilde{h}_j||_{L^2} \to 0$. Consequently, 
$\Pi_\e \tilde{h}_j \to h_j \in K_j$, at least uniformly on compact sets of $M_j$, and 
$\Pi_\e \tilde{h}_j \to 0$ on the other component $M_{j'}$, $j' \neq j$. In particular, 
the set of all $\Pi_\e \tilde{h}_j$, as $h_j$ varies over $K_j$, is a basis of $V_\e$.

We define the space $\rho^\mu \Lambda^{k,\al}_0 S^2_0(M_\e)^\perp$  to be the closed subspace of 
elements of $\rho^\mu \Lambda^{k,\al}_0 S^2_0(M_\e)$ which are $L^2$ orthogonal to elements of 
$V_\e$. We prove 
\begin{proposition} 
For $\e$ sufficiently small and $0 < \mu < \inf (2, n/2)$, the mapping
\[
\hat L_{g_\e} : \calS \oplus \rho^\mu \Lambda^{2,\al}_0 S^2_0(M_\e)^\perp \longrightarrow 
\rho^\mu \Lambda^{0,\al}_0 S^2_0(M_\e) 
\]
is an isomorphism and its inverse $\hat G_{g_\e}$ is bounded  independently of $\e$.
\label{pr:8}
\end{proposition}
{\bf Proof:} Following the proof of Proposition~\ref{pr:7} we deduce that 
\[
L_{g_\e} : \rho^\mu \Lambda^{2,\al}_0 S^2_0(M_\e)^\perp  \longrightarrow   
\rho^\mu \Lambda^{0,\al}_0 S^2_0(M_\e)^\perp 
\]
is an isomorphism whose inverse is bounded independently of $\e$. 

Next we show that for $\e$ small, and all $h \in V_\e$, there exists $\gamma \in \calS$ for which 
\[
\int_{M_\e} \langle \hat L_{g_\e} (\gamma, 0) , h \rangle \neq  0
\]
This follows from the fact that for any $\gamma \in \calS $, $\hat L_{g_\e} (\gamma, 0)
 = \hat L_{g_1} (\gamma^{(1)}, 0) + \hat L_{g_2} (\gamma^{(2)}, 0)$, and the first term on 
the right is supported in $M_1 - \{p_1\}$ and the second on $M_2 - \{p_2\}$. 
Thus we need to find $\gamma^{(1)}$ and $\gamma^{(2)}$ such that
\[
\int_{M_\e} \left( \langle \hat L_{g_1} (\gamma^{(1)}, 0), h \rangle + \langle \hat L_{g_2} 
(\gamma^{(2)},0) ,h  \rangle \right)\neq 0.
\]
It is enough to check this for all $h_\e = \Pi_\e \tilde{h}$, where $\tilde h$ is the 
extension to all of $M_\e$ of an arbitrary element of $K_1$ or $K_2$. 
Letting $\e \to 0$ on the left gives
\begin{equation}
\int_{M_1}\langle\hat L_{g_1}(\gamma^{(1)},0), h_1 \rangle + \int_{M_1} 
\langle \hat L_{g_2} (\gamma^{(2)},0) ,h_2  \rangle 
\label{eq:gggg}
\end{equation}
for some elements $h_j \in K_j$, not both equal to $0$. That we can choose 
$\gamma^{(j)}$ so that this can be made nonvanishing is the content of the 
modification of Proposition~\ref{pr:qqq} at the end of \S 2. The proof is complete. \hfill $\Box$

\section{Proof of Theorem 1}

Following the development of the linear analysis in the last section, it is now a simple matter to complete the proof of the main theorem. Recall that our goal is to find a correction term $k_\e$ to $g_\e$ so that $g_\e + k_\e$ is Poincar\'e-Einstein. The nonlinear operator $k  \rightarrow  N_{g_\e}(k)$ is a second order quasilinear operator with coefficients which are polynomial in the entries of $(\nabla^{g_\e})^j \, k$, $j=0,1,2$, and $(g_\e +k)^{-1}$. The same is true for the remainder term
\[
Q_{g_\e} (k) : =  N_{g_\e}(k) - N_{g_\e}(0) - L_{g_\e} k,
\]
and in addition, if $\mu >0$, there exists a constant $C>0$, which does not depend on $\e$, such that 
\begin{equation}
|| Q_{g_\e}(k_2) - Q_{g_\e}(k_1)||_{0, \alpha, \mu} \leq C \, \left( ||k_2||_{2, \alpha, \mu} + ||k_1||_{2, \alpha, \mu} \right) \, ||k_2 -k_1||_{2, \alpha, \mu},
\label{eq:mlml}
\end{equation}
for all $k_1, k_2 \in \rho^\mu \, \Lambda^{2, \alpha}_0 \, S^2_0 (M_\e)$ satisfying $||k_1||_{2, \alpha, \mu} + ||k_2||_{2, \alpha, \mu} \leq 1$. 

In the case where the metrics $g_1$ and $g_2$ are both nondegenerate, we fix $\mu \in (0, 1)$ and use the result of Proposition~\ref{pr:6} to shows that 
\begin{equation}
|| N_{g_\e} (0) ||_{0, \alpha, \mu}\leq C'\, \e^{2-\mu}
\label{eq:eee}
\end{equation}
Then, Proposition~\ref{pr:7} allows use to rephrase the equation $N_{g_\e} (k) =0$  as a fixed point problem
\[
k = - \, G_{g_\e} \, \left(N_{g_\e}(0) + Q_{g_\e} (k) \right).
\]
The fact that $G_{g_\e}$ is uniformly bounded together with (\ref{eq:mlml}) and (\ref{eq:eee}) implies that the mapping $k  \rightarrow  - \, G_{g_\e} \, \left(N_{g_\e}(0) + Q_{g_\e} (k) \right)$ is a contraction mapping from a small ball of $\rho^\mu \, \Lambda^{2, \alpha}_0 \, S^2_0 (M_\e)$ into itself. This complete the proof of the existence of a solution of $N_{g_\e} (k) =0$.

The case where one of the metrics $g_1$ or $g_2$ is degenerate can be treated similarly using 
Proposition~\ref{pr:8} instead of Proposition~\ref{pr:7}. Observe that an estimate similar to 
(\ref{eq:mlml}) is valid for all $(\gamma_1,k_1), (\gamma_2, k_2) \in {\calS} \oplus \rho^\mu \, 
\Lambda^{2, \alpha}_0 \,  S^2_0 (M_\e)$. We leave the details to the reader.

\section{Scalar positive conformal infinities} 

Given any two \PE metrics $(M_1,g_1)$ and $(M_2,g_2)$, we have shown  how to produce a family of \PE metrics on the boundary connected sum  $M_1 \#_b M_2$ in such a way that the new conformal infinity  is  unchanged away from a small neighborhood of the gluing points,  in case the linear Einstein operators on both factors have trivial  $L^2$ kernels, or is altered only very slightly, in the general case.  {}From here it is only a small step to prove Corollary~1 concerning how  the Yamabe class of the conformal infinities of these new metrics  relate to the initial conformal infinities $\frakc(g_1)$ and   $\frakc(g_2)$. 

Let us change notation slightly from the last section and write the \PE metric as $g_\e$, and decompose
it as a sum of an explicit  approximate solution $\tg_\e$ and the deformation term $k_\e$. Let us 
denote by $h_1$, $h_2$ and $h_\e = \thh_\e + q_\e$ the explicit metric representatives of the conformal 
infinities of each of these metrics, where $\thh_\e$ is the conformal infinity of $\tg_\e$ and, 
in the  degenerate case, $q_\e \in \calS$, but equals $0$ otherwise. For simplicity, we also write  
$Y_j = \del M_j$.

We first recall the connected sum construction for metrics of constant scalar curvature from 
\cite{MPU} and \cite{J}. The idea of the  construction and much of the implementation is almost 
exactly the same  as what we have done here; indeed, the main substantial difference is  the need 
for the theory of uniformly degenerate operators for the  interior problem. In any case, suppose 
that $h_1$ and $h_2$ are  metrics of constant scalar curvature on $Y_1$ and $Y_2$, respectively; 
\cite{MPU} handles the case where the scalar curvatures are the same  positive number, while 
\cite{J} treats the more general situation  where the constants may differ and possibly even have 
different signs.  Having chosen points $p_j \in Y_j$, one identifies by inversion the  small annuli
of inner and outer radii $\e$ and $2\e$,  for example,  around these points to define $Y_1 \# Y_2$, 
and then uses a  partition of  unity to patch together the metrics to define a family  of metrics 
$\thh_\e$. (This construction is phrased differently in \cite{MPU}: there, fixed annuli around the
$p_j$ are transformed conformally to  long cylinders, and these are then patched together; the 
approximate  solution metric $\thh_\e$ is given by a conformal factor which has the shape of a 
$\cosh$ curve on this cylinder, hence is exponentially in the middle small relative to its length. 
The equivalence of this picture with the other one is immediate.)  This step is clearly identical 
to what is happening on the boundary in our construction of the approximate \PE metrics $\tg_\e$. 

The constant scalar curvature metric $h_\e$ is obtained by solving the scalar Yamabe equation, and is a conformal deformation from the background metric $\thh_\e$. In particular, if $h_1$ and $h_2$ are  both scalar positive metrics, then the conformal class of $\thh_\e$ on  $Y_1 \# Y_2$ is also scalar positive, provided $\e$ is small enough.  

In the case of nondegenerate \PE gluing, the conformal infinity of $g_\e$ is the same as that of $\tg_\e$, and we have just shown that  this is scalar positive if this is true for both summands. In the  degenerate case, the proof of Theorem 1 shows that $\frakc(g_\e)$ is a  $\calC^2$ small perturbation of $\frakc(\tg_\e)$, and since the Yamabe  functional is continuous in the $\calC^2$ topology, the conformal  class $\frakc(g_\e)$ is again scalar positive. This concludes the  proof of Corollary 1.

\section{Plumbing and surgery?}

We have proved that it is possible to perform a boundary connected sum  in the category of 
\PE metrics. There are many other interesting ways  to join together $M_1$ and $M_2$ along 
their boundaries as differentiable manifolds, and it is natural to ask whether these operations 
may also be done in the \PE category.  The two operations we have in mind are: 

\medskip

$\bullet$ Boundary plumbing~: Suppose that $\Sigma$ is a $k$-dimensional  manifold which is smoothly 
embedded in $\del M_j$, $j=1,2$, in such a  way that the normal bundles $\nu_j$ of these embeddings  
are equivalent. If $S\nu_j$ denote the unit sphere bundles, then the  identification $S\nu_1 
\cong S\nu_2$ extends to an orientation reversing bundle map  $\nu_1\setminus \{0\} \to 
\nu_2 \setminus \{0\}$ which 
is homogeneous of  degree $-1$ on the fibres, which we call the inversion $I$. The bundles $\nu_j 
\oplus \RR^+$ are diffeomorphic to inward-pointing tubular  neighborhoods $\calT_j^+$ of $\Sigma$ 
in $M_j$ (both are half-ball  bundles over $\Sigma$). The natural extension of $I$ is a  
diffeomorphism $\calT_1^+ - \Sigma \to \calT_2^+ - \Sigma$, and using  this we define the 
boundary join of $M_1$ and $M_2$ along $\Sigma$,  $M_1 \#_{b,\Sigma} M_2$.

\medskip

$\bullet$ Boundary surgery~: Suppose that $S^k \subset \del M_1$ and  $S^{n-k-1} \subset  \del M_2$ are 
spheres, both with trivial normal  bundles, and let $\calT_j$, $\calT_j^+$, be the corresponding 
tubular  neighborhoods in $\del M_j$, $M_j$, respectively.  Thus  $M_j - \calT_j^+$ are both 
manifolds with corners, each with two  boundary  hypersurfaces
\[
(\del M_1 - \calT_1) \cup (S^k \times B^{n-k})  \qquad \mbox{and} \qquad (\del M_2 - \calT_2) 
\cup (B^{k+1}\times S^{n-k-1}),
\]
respectively. A standard construction in topology joins these two spaces using the plug 
$B^{k+1}\times B^{n-k}$, which has boundary $(S^k \times B^{n-k}) \cup (B^{k+1}\times S^{n-k-1})$, 
and hence may be inserted between the two summands $M_j - \calT_j^+$ to define the surgered manifold 
$M_1 \#_{b,\sigma} M_2$ (the $b,\sigma$ subscript simply means  `boundary surgery', but we 
suppress the dimension of the surgery from the notation). 

\medskip

Let $(M_j,g_j)$ be two \PE metrics.  We ask the following questions:
\begin{itemize}
\item Suppose $\Sigma \subset \del M_j$ and the two normal bundles  $\nu_j$ are equivalent. 
Does the boundary join  $M_1 \#_{b,\Sigma} M_2$ admit a family of \PE metrics which converges 
nicely to $g_j$ on compact sets of $\overline{M}_j - \Sigma$, and such  that the conformal
infinity $\frakc(g_\e)$ is close to $\frakc(g_j)$  away from the neck region~? 
\item If $S^k$ and $S^{n-k}$ are spheres with trivial normal bundles  in $\del M_1$ and 
$\del M_2$, respectively, then does the surgered  manifold $M_1 \#_{b,\sigma} M_2$ admit a 
similar family of \PE metrics  $g_\e$~?
\end{itemize}

The utility of these constructions is obvious, and in particular if the second were always possible, 
then it would not be unreasonable to hope that any (compact, nullcobordant, simply connected) 
scalar-positive manifold might admit a \PE filling. There would be many other applications too. 
We are not able to answer either of these questions one way or the other, but these are clearly 
important directions for future research. However, we suspect that if either $M_1 \#_{b,\Sigma} M_2$ 
or $M_1 \#_{b,\sigma} M_2$ admits a \PE metric  $g$, then $g$ is `quite far' from any family of 
`locally constructed' approximate solution metrics $\tilde{g}_\e$, unlike the construction above.
In particular, it does not seem likely that there should exist families of \PE metrics $g_\e$ in 
either case which converge to $g_1$  and $g_2$ away from the necks and which have restrictions to 
the necks which are close to any simple model form. As heuristic evidence for this, we note
that were such a family $g_\e$ to exist, then one would expect that rescalings of its
restriction to the neck region should converge to some sort of model Einstein metric.
For example, in the boundary connected sum construction, this model metric is just hyperbolic space.
However, calculations seem to indicate that there there are no warped product candidates for the
model Einstein metrics in these more general cases, and it is not clear where else to look.
An additional nonrigorous counterargument is that if $g_\e$ were to have negative sectional
curvature in the neck region, then one would be able to join or surger together two
copies of hyperbolic space and obtain a manifold with negative sectional curvature
everywhere, but which is simply connected and has nontrivial higher homotopy groups, which
is impossible.

\end{document}